\documentclass[sn-mathphys-num]{sn-jnl}

\usepackage{graphicx}%
\usepackage{multirow}%
\usepackage{amsmath,amssymb,amsfonts}%
\usepackage{amsthm}%
\usepackage{mathrsfs}%
\usepackage[title]{appendix}%
\usepackage{xcolor}%
\usepackage{textcomp}%
\usepackage{manyfoot}%
\usepackage{booktabs}%
\usepackage{algorithm}%
\usepackage{algorithmicx}%
\usepackage{algpseudocode}%
\usepackage{listings}%

\theoremstyle{thmstyleone}%
\newtheorem{theorem}{Theorem}
\theoremstyle{thmstyletwo}%
\newtheorem{remark}{Remark}%

\theoremstyle{thmstylethree}%
\newtheorem{definition}{Definition}%

\begin{document}

\title[Pandey-Upadhyay's wavelet transform]
{Pandey-Upadhyay's wavelet transform and 
microlocal Sobolev singularities of functions}


\author[1]{\fnm{Akira} \sur{Lee}}\email{waa\_ri@cc.nara-wu.ac.jp}

\author*[2]{\fnm{Shinya} \sur{Moritoh}}\email{moritoh@cc.nara-wu.ac.jp}

\affil[1]{Graduate School of Humanities and Sciences, 
Nara Women's University, 
Nara, Japan 630-8506}

\affil*[2]{Department of Mathematics, Nara Women's University, 
Nara, Japan 630-8506~~ (Tel: +81-742-20-3981)}



\abstract{
The aim of the paper is to define the microlocal Sobolev 
singularities of functions using Pandey-Upadhyay's 
wavelet transform and provide a comparison with H\"ormander's 
microlocal 
singularities.}


\maketitle

\vspace{-1cm}

\hspace{0.7cm}

\section{Introduction}

The aim of the paper is to define the microlocal Sobolev 
singularities of functions using Pandey-Upadhyay's 
wavelet transform and provide a comparison 
with H\"ormander's 
microlocal singularities.
 
Let $f$ be a square integrable function, i.e. $f\in L^2({\mathbb R}^n)$. 
Let the wavelet function $\psi$ be rapidly decreasing, i.e. $\psi
 \in {\mathcal S}({\mathbb R}^n)$, and its Fourier transform 
$\hat{\psi}$ be compactly supported in the cube $[1/2,2]^n$. 
Then the wavelet transform of $f$ is defined as follows (see 
Pandey-Upadhyay \cite{PU}): 
\vspace{-1mm}
\[
\begin{array}{l}
\displaystyle
W_{\psi}f(x,\xi)
=\int_{{\mathbb R}^n}
f(t)\, |\xi_1 \cdots \xi_n| \overline{\psi}
(\xi_1(t_1-x_1),\ldots,\xi_n(t_n-x_n))dt, 
\end{array}
\]
where $x=(x_1,x_2,\ldots,x_n) \in{\mathbb R}^n, \,
\xi=(\xi_1,\xi_2,\ldots,\xi_n)\in({\mathbb R}-\{0\})^n$. 
The inversion formula for Pandey-Upadhyay's 
wavelet transform reads as follows: 
\begin{equation}\label{30}
\begin{array}{l}
\displaystyle
f(t)
=C_{\psi}^{-1}\iint_{{\mathbb R}^n\times({\mathbb R}-\{0\})^n}
\!\!\!
W_{\psi}f(x,\xi) |\xi_1 \cdots \xi_n| \times
\vspace{0.2cm}
\\
\displaystyle
\hspace{3.5cm}
\times \psi
(\xi_1(t_1-x_1), \ldots, \xi_n(t_n-x_n))dxd\xi/|\xi_1\cdots\xi_n|, 
\end{array}
\end{equation}
where $C_{\psi}=
\int_{{\mathbb R}^n}|\hat{\psi}(\omega)|^2
d\omega/|\omega_1\cdots\omega_n|$ 
is a finite, nonzero constant. 

It is easily seen that the Fourier transform of 
$W_{\psi}f(x,\xi)$ with respect to $x$ can be calculated as follows: 
\[
\left(W_{\psi}f(\cdot,\xi)\right)^{\wedge}(\tau):=
\int_{{\mathbb R}^n}W_{\psi}f(x,\xi)e^{-ix\cdot \tau}dx
= 
\hat{f}(\tau)\overline{\hat{\psi}}
(\xi_1^{-1}\tau_1,\ldots, \xi_n^{-1}\tau_n). 
\]
Since $\hat{\psi}$ is compactly supported in the cube $[1/2,2]^n$, 
the function 
$\left(W_{\psi}f(\cdot,\xi)\right)^{\wedge}(\tau)$ is compactly 
supported in rectangular parallelepiped 
$[(1/2)\xi_1,2\xi_1]\times\cdots\times[(1/2)\xi_n,2\xi_n]$. 

We note that, in the definitions of wavelet transforms, 
Pandey-Upadhyay \cite{PU} uses 1-dimensional dilations 
in all coordinate directions, 
while Moritoh \cite{M} uses $n$-dimensional rotation 
and 1-dimensional dilation in the radial direction.
\vspace{1mm}

In the following section 2, as an anlogue of the 
singularities considered in 
Moritoh \cite{M}, we define the microlocal Sobolev 
singularities of functions using Pandey-Upadhyay's 
wavelet transform, and state our main theorem 
providing a comparison with H\"ormander's 
microlocal singularities. In the proof, we 
will focus on the corresponding key parts of \cite{M}.

\section{Definitions, main theorem, and the proof}

Let us first recall the definition of H\"ormander's microlocal 
singularities. 
\begin{definition}[H\"ormander \cite{H}]
Let $f\in L^2({\mathbb R}^n)$, 
$(x_0, \xi^0)\in {\mathbb R}^n\times({\mathbb R}^n-\{0\})$, 
and $s\ge0$. 
Then, $f$ is said to belong to the  
Sobolev space 
$H^s$ microlocally at $(x_0,\xi^0)$ in the sense of 
H\"ormander if
there exist a cutoff function 
$\phi\in C_0^{\infty}({\mathbb R}^n)$ identically 
equal to 1 in a neighborhood of $x_0$ and a conical 
neiborhood $\Gamma(\xi^0)$ of $\xi^0$ such that 
\begin{equation}\label{2}
\int_{\Gamma(\xi^0)}|(\phi f)^{\wedge}(\xi)|^2(1+|\xi|^2)^s d\xi<\infty.
\end{equation}
\end{definition}

\noindent
We now define the microlocal singularities of functions 
using Pandey-Upadhyay's wavelet transform.
\begin{definition}
Let $f\in L^2({\mathbb R}^n)$, 
$(x_0, \xi^0)\in {\mathbb R}^n\times ({\mathbb R}-\{0\})^n$, 
and $s\ge0$. Then, $f$ is said to belong to the  
Sobolev space 
$H^s$ microlocally at $(x_0,\xi^0)$ in the 
sense of PU if
there exist a neiborhood $U(x_0)$ of $x_0$ and a conical 
neiborhood $\Gamma(\xi^0)$ of $\xi^0$ such that 
\begin{equation}\label{1}
\iint_{U(x_0)\times\Gamma(\xi^0)}|W_{\psi}f(x,\xi)|^2(1+|\xi|^2)^sdxd\xi/ |\xi_1 \cdots \xi_n|<\infty. 
\end{equation}
\end{definition}

\vspace{1mm}
In order to state the main theorem, let us assume, 
for simplicity, that $n=2$ and that 
both coordinates of $\xi^0=(\xi^0_1, \xi^0_2)$ are positive, i.e. 
$\xi^0\in({\mathbb R}_{>0})^2$, and 
define the conical neighborhood 
$\Gamma_{\psi}(\xi^0)$ of $\xi^0$, with some margin 
depending on the wavelet $\psi$, as follows.  
\begin{equation}\label{3}
\Gamma_{\psi}(\xi^0)=\{\tau=(\tau_1,\tau_2)\in({\mathbb R}_{>0})^2|\,
(1/5)(\xi^0_1/\xi^0_2) <\tau_1/\tau_2< 5(\xi^0_1/\xi^0_2)\}. 
\end{equation}
\vspace{1mm}

\begin{theorem}
Let $f\in L^2({\mathbb R}^2)$, $(x_0, \xi^0)\in {\mathbb R}^2\times 
({\mathbb R}_{>0})^2$, and $s\ge0$.
\vspace{1mm}

\noindent
{\rm (a)} 
If {\rm (\ref{2})} holds for a cutoff function 
$\phi\in C_0^{\infty}({\mathbb R}^2)$ identically 
equal to 1 in a neighborhood of $x_0$ and the conical neighborhood 
$\Gamma_{\psi}(\xi^0)$ of $\xi^0$ defined by {\rm (\ref{3})}, 
then $f$ belongs to the  
Sobolev space 
$H^s$ microlocally at $(x_0,\xi^0)$ in the sense of PU. 
\vspace{1mm}

\noindent
{\rm (b)}
Conversely, if {\rm (\ref{1})} holds for a neighborhood 
$U(x_0)$ of $x_0$ and the conical neighborhood 
$\Gamma_{\psi}(\xi^0)$ 
defined by {\rm (\ref{3})}, then $f$ belongs to the Sobolev space 
$H^s$ microlocally at $(x_0,\xi^0)$ in the sense of H\"ormander. 
\end{theorem}
\vspace{1mm}

\begin{remark}
In Moritoh \cite{M}, a comparison theorem is stated 
using wave front sets and considers enlarging the 
wave front sets only in frequency space. 
This corresponds to how $\Gamma_{\psi}(\xi^0)$ is defined by 
(\ref{3}). We also note that $\Gamma_{\psi}(\xi^0)$ needs to 
be a conical neighborhood that is slightly larger than 
the smallest conical one containing 
$[(1/2)\xi_1^0,2\xi_1^0]\times[(1/2)\xi_2^0, 2\xi_2^0]$.
\end{remark}

\begin{proof}
We use the method employed in the proof of 
theorem 1 in \cite{M}. 
We also note that proposition 2 (local property) 
in \cite{M} also holds true for Pandey-Upadhyay's 
wavelet transform: 
If $x_0$ does not 
belong to $\text{supp}\,f$, then there exists a neighbourhood 
$U(x_0)$ of $x_0$ such that $W_{\psi}f(x,\xi)$ is rapidly 
decreasing in $\xi$ uniformly in $x\in U(x_0)$.

First, let us prove (a). Our assumption is that, for a cutoff function 
$\phi\in C_0^{\infty}({\mathbb R}^2)$ identically 
equal to 1 in a neighborhood of $x_0$, the following holds true.  
\begin{equation}\label{10}
\int_{\Gamma_{\psi}(\xi^0)}|(\phi f)^{\wedge}(\xi)|^2(1+|\xi|^2)^s d\xi<\infty.
\end{equation}
Let a neighborhood $U(x_0)$ of $x_0$ be contained in the set 
$\{x\in{\mathbb R}^2| \phi(x)=1\}$ and let 
a conical neighborhood $\Gamma(\xi^0)$ of $\xi^0$ be 
defined by 
$\{\tau=(\tau_1,\tau_2)\in({\mathbb R}_{>0})^2|\,
(4/5)(\xi^0_1/\xi^0_2) <\tau_1/\tau_2< (5/4)(\xi^0_1/\xi^0_2)\}$. 
Using the same argument as for the first half of 
the proof of theorem 1 in \cite{M}, it suffices to 
show only the following.
\[
\iint_{{\mathbb R}^2\times\Gamma(\xi^0)}
|W_{\psi}(\phi f)(x,\xi)|^2(1+|\xi|^2)^s dxd\xi/ |\xi_1\xi_2| 
<\infty.
\]
This can be rewritten as 
\begin{equation}\label{ineq1}
\iint_{{\mathbb R}^2\times\Gamma(\xi^0)} 
|{(\phi f)}^{\wedge}(\tau) \hat{\psi}(\xi_1^{-1}\tau_1, 
\xi_2^{-1}\tau_2)|^2(1+|\xi|^2)^sd\tau d\xi/ |\xi_1\xi_2|
<\infty.
\end{equation}
Changing the variable $\xi$ to the variable $\omega$, where  
$\omega_1=\xi_1^{-1}\tau_1$ and $\omega_2=\xi_2^{-1}\tau_2$, 
and taking the support conditon for $\hat{\psi}(\omega)$ into 
account, we can see that $d\xi/ |\xi_1\xi_2|=
d\omega/ |\omega_1\omega_2|$ and $1+|\xi|^2 \le 4(1+|\tau|^2)$. 
Moreover,  when $\xi\in\Gamma(\xi^0)$, from the equality 
$\tau_1/\tau_2=(\omega_1/\omega_2)(\xi_1/\xi_2)$ we see that 
$\tau_1/\tau_2$ is between (1/4)(4/5) and 4(5/4) times the 
ratio $\xi^0_1/\xi^0_2$, i.e. 
$\tau\in \Gamma_{\psi}(\xi^0)$. 
Therefore, the left side of (\ref{ineq1}) is bounded from above by 
\[
\int_{\Gamma_{\psi}(\xi^0)} 
|{(\phi f)}^{\wedge}(\tau)|^2\,4^s (1+|\tau|^2)^s 
\left(\int_{{\mathbb R}^2}|\hat{\psi}(\omega_1, \omega_2)|^2
d\omega/|\omega_1\omega_2|\right) d\tau, 
\]
which is finite from our assumption (\ref{10}).

Next, let us prove (b). Our assumption is that, for 
a neighborhood $U(x_0)$ of $x_0$, the following holds true.  
\begin{equation}\label{20}
\iint_{U(x_0)\times\Gamma_{\psi}(\xi^0)}
|W_{\psi}f(x,\xi)|^2(1+|\xi|^2)^sdxd\xi/ |\xi_1\xi_2|<\infty. 
\end{equation}
Let $\phi\in C_0^{\infty}({\mathbb R}^2)$ 
be compactly supported in $U(x_0)$ and be identically equal to 1 
in some neighborhood of $x_0$.  Then, taking the inversion formula 
(\ref{30}) into account, and  
using the same argument as for the second half of 
the proof of theorem 1 in \cite{M}, it suffices to 
show only the following:  The function defined as 
\begin{equation}\label{4}
\begin{array}{l}
f_1(t)
\\
\displaystyle
=C_{\psi}^{-1}\!\!\! \iint_{{\mathbb R}^2\times\Gamma_\psi(\xi^0)}
\!\!\!\!\!\!\!\!\!
\phi (x)W_\psi f(x,\xi) |\xi_1\xi_2| \psi
(\xi_1(t_1-x_1), \xi_2(t_2-x_2))dxd\xi/|\xi_1\xi_2|
\end{array}
\end{equation}
belongs to the Sobolev space $H^s({\mathbb R}^2)$. 
 
By putting 
$g(x,\xi)=\phi (x)W_\psi f(x,\xi)(1+|\xi|^2)^{s/2}$, 
our assumption (\ref{20}) can be rewritten as 
\begin{equation}\label{5}
\iint_{{\mathbb R}^2\times\Gamma_\psi(\xi^0)}
|g(x,\xi)|^2 dxd\xi/|\xi_1\xi_2| < \infty.
\end{equation}
The Fourier transform of (\ref{4}) is equal to 
\[
\widehat{f_1}(\tau)=C_{\psi}^{-1}\int_{\Gamma_\psi(\xi^0)}
\left(g(\cdot,\xi)\right)^{\wedge}(\tau)
\hat{\psi}(\xi_1^{-1}\tau_1, \xi_2^{-1}\tau_2)
(1+|\xi|^2)^{-s/2}d\xi/|\xi_1\xi_2|. 
\]
Therefore, by using the Schwarz inequality, the same 
argument immediately following (\ref{ineq1}), 
the Fubini theorem, the Plancherel theorem and 
our assumption (\ref{5}), we obtain  

\[
\begin{array}{l}
\displaystyle
\int_{{\mathbb R}^2}|\widehat{f_1}(\tau)|^2(1+|\tau|^2)^sd\tau
\vspace{0.2cm}
\\
\displaystyle
\le 
C_{\psi}^{-2}\int_{{\mathbb R}^2}
\left(\int_{\Gamma_\psi(\xi^0)}
|\left(g(\cdot,\xi)\right)^{\wedge}(\tau)|^2 
d\xi/|\xi_1\xi_2| \right) \times
\\
\displaystyle
\hspace{3.5cm}
\times \left(\int_{{\mathbb R}^2}
|\hat{\psi}(\xi_1^{-1}\tau_1, \xi_2^{-1}\tau_2)|^2
\left(\frac{1+|\tau|^2}{1+|\xi|^2}\right)^s 
d\xi/|\xi_1\xi_2| \right)
d\tau
\vspace{0.2cm}
\end{array}
\]
\[
\begin{array}{ll}
\displaystyle
\le
C_{\psi}^{-2}
\int_{{\mathbb R}^2}
\left(\int_{\Gamma_\psi(\xi^0)}
|\left(g(\cdot,\xi)\right)^{\wedge}(\tau)|^2
d\xi/|\xi_1\xi_2|\right)\times
\\
\displaystyle
\hspace{3.5cm}
\times 
4^s
\left(\int_{{\mathbb R}^2}|\hat{\psi}(\omega_1, \omega_2)|^2
d\omega/|\omega_1\omega_2|\right)
d\tau
\vspace{0.2cm}
\\
\displaystyle
=
4^sC_{\psi}^{-1}
\int_{\Gamma_\psi(\xi^0)}
\left(\int_{{\mathbb R}^2}
|g(x,\xi)|^2
dx\right)
d\xi/|\xi_1\xi_2|<\infty. 
\end{array}
\]
\end{proof}

\vspace{1mm}

\begin{remark}
From the above, it is easy to see that 
function spaces with dominating mixed smoothness are 
worth considering by using Pandey-Upadhyay's 
wavelet transform. 
\end{remark}

\begin{remark} 
We note the paper \cite{HPS} in relation to our study. 
In Han-Paul-Shukla \cite{HPS}, 
continuous shearlet transforms are used 
to characterize the microlocal Sobolev wave front sets, 
the 2-microlocal spaces, and local Hölder spaces.  See also 
the original reference Guo-Kutyniok-Labate \cite{GKL}.
\vspace{4mm}
\end{remark}

\section*{Acknowledgements}
We would like to express our sincere gratitude 
to Professor Kensuke Fujinoki (Kanagawa University) 
for his helpful comments and suggestions on this 
paper at the Wavelet Seminar organized by him.

\vspace{2mm}





\end{document}